\documentclass[12pt]{article} 
\usepackage{fourier}
\usepackage{index}
\usepackage[expansion=false]{microtype}
\usepackage[reqno]{amsmath}
\usepackage{fixmath}
\usepackage{amssymb, amsthm, enumerate}
\usepackage{mathtools}

\usepackage{hyperref}

\newtheoremstyle{plainsl}%
	{\topsep}
	{\topsep}
	{\slshape} 
	{}
	{\normalfont\bfseries}
	{.}
	{ }
	{}

\swapnumbers

{\theoremstyle{plainsl}
\newtheorem{theorem}{Theorem}[section]
\newtheorem{lemma}[theorem]{Lemma}
}
{\theoremstyle{remark}
}

\newcommand\lref[1]{Lemma~\ref{lem:#1}}

\newcommand\cref[1]{Corollary~\ref{cor:#1}}

\renewcommand\proof{\noindent\textsl{Proof. }}
\newcommand\sqr[2]{{\vbox{\hrule height.#2pt
    \hbox{\vrule width.#2pt height#1pt \kern#1pt
        \vrule width.#2pt}\hrule height.#2pt}}}
\renewcommand\qed{%
	\ifmmode\eqno\sqr53
	\else\nolinebreak\ \hfill\sqr53\medbreak\fi}

\numberwithin{equation}{section}

\newcommand\al{\alpha}
\newcommand\be{\beta}
\newcommand\de{\delta}
\newcommand\De{\Delta}

\newcommand\ga{\gamma}

\newcommand\sg{\sigma}

\renewcommand\th{\theta} 

\newcommand\cA{{\mathcal A}}

\newcommand\cx{{\mathbb C}}
\newcommand\fld{{\mathbb F}}

\newcommand\ints{{\mathbb Z}}
\newcommand\re{{\mathbb R}}

\newcommand\comp[1]{{\mkern2mu\overline{\mkern-2mu#1}}}
\newcommand\diff{\mathbin{\mkern-1.5mu\setminus\mkern-1.5mu}}

\newcommand\seq[3]{#1_{#2},\ldots,#1_{#3}}

\newcommand\pmat[1]{\begin{pmatrix} #1 \end{pmatrix}}
\newcommand\sm[3]{\sum_{#1=#2}^{#3}}

\newcommand\one{{\bf1}}

\newcommand\ip[2]{\langle #1,#2\rangle}
\newcommand\hm{\widehat M}
\newcommand\hmx{\hm_X}
\newcommand\psd{\succcurlyeq}

\title{Average Mixing of Continuous Quantum Walks} 

\author{
	Chris Godsil\\
	Combinatorics \& Optimization\\
	University of Waterloo}

\begin{document}
\maketitle
	
\begin{abstract}
    If $X$ is a graph with adjacency matrix $A$, then we define $H(t)$ to be
    the operator $\exp(itA)$. The Schur (or entrywise) product $H(t)\circ H(-t)$
    is a doubly stochastic matrix and because of work related to quantum computing,
    we are concerned with the \textsl{average mixing matrix} $\hmx$, defined by
    \[
        \hmx = \lim_{C\to\infty}\frac1C \int_0^C H(t)\circ H(-t) \,dt.
    \]
    In this paper we establish some of the basic properties
    of this matrix, showing that it is positive semidefinite and that its entries 
    are always rational. We see that in a number of cases its form is surprisingly simple.
    Thus for the path on $n$ vertices it is equal to
    \[
        \frac1{2n+2}(2J+I+T)
    \]
    where $T$ is the permutation matrix that swaps $j$ and $n+1-j$ for each $j$.
    If $X$ is an odd cycle or, more generally, if $X$ is one of the graphs in 
    a pseudocyclic association scheme on $n$ vertices with $d$ classes, each of
	valency $m$, then its average mixing matrix is
    \[
        \frac{n-m+1}{n^2}J + \frac{m-1}{n}I.
    \]
    (One reason this is interesting is that a graph in a pseudocyclic scheme may
    have trivial automorphism group, and then the mixing matrix is more symmetric than
	the graph itself.)
\end{abstract}

\section{Average Mixing}

Let $X$ be a graph with adjacency matrix $A$. We define a transition matrix $H_X(t)$ 
by
\[
    H_X(t) := \exp(itA).
\]
For a physicist this matrix determines a continuous quantum walk. We note that it
is both symmetric and unitary, in particular $\comp{H_X(t)}=H_X(-t)$.
For relevant recent surveys see, e.g., Kendon and Tamon \cite{KendTam2010}, 
Godsil \cite{Godsil2011a}.

The matrix $H(t)$ gives rise to a family of probability densities as follows.
Let $A\circ B$ denote the \textsl{Schur product} of two matrices $A$ and $B$
with the same order. Thus
\[
    (A\circ B)_{u,v} = A_{u,v}B_{u,v}.
\]
and we will use $A^{\circ 2}$ to denote $A\circ A$. Then if
\[ 
    M_X(t) := H_X(t)\circ H_X(-t)
\]
we see that $M_X(t)$ is a nonnegative real matrix and, since $H_X(t)$ is unitary,
each row and column of $M_X(t)$ sums to 1. We use $e_u$ to denote the
standard basis vector of $\cx^{V(X)}$ indexed by the vertex $u$, thus
we have a family of probability densities $e_u^TM_X(t)$, and we are concerned with
the behavior of these densities.
In this paper we are interested in the matrix $\hmx$, which we define by
\[
	\hmx = \frac1C \int_0^C M_X(t)\, dt.
\]
Following \cite{Adamczak2008}, we call this the \textsl{average mixing matrix} of $X$. 

To work with $M_X$ and $\hmx$, we use the spectral decomposition of $A$. This allows
us to write $A$ as
\[
    A = \sum_r \th_r E_r
\]
where $\th_r$ runs over the distinct eigenvalues of $A$ and $E_r$ is the matrix
representing orthogonal projection onto the eigenspace belonging to $\th_r$.
Then we have
\[
    H_X(t) = \sum_r \exp(i\th_r t)E_r
\]
and 
\[
    M_X(t) = \sum_r E_r^{\circ 2} +2\sum_{r<s}\cos((\th_r-\th_s)t)\, E_r\circ E_s.
\]
From this we have:

\begin{lemma}
    If $A=\sum_r \th_r E_r$ is the spectral decomposition of $S=A(X)$, then
    \[
        \hmx = \sum_r E_r^{\circ2}.\qed
    \]
\end{lemma}

\section{Properties of the Average Mixing Matrix}

In this section we derive some basis properties of the average mixing matrix.
We know already that it is symmetric. By a famous theorem
of Schur, the Schur product $M\circ N$ of two positive semidefinite matrices
is positive semidefinite and, since the sum of positive semidefinite matrices
is positive semidefinite, we see that $\hmx$ is positive semidefinite.
If $A$ and $B$ are symmetric matrices of the same order, we write $A\psd B$
to denote that $A-B$ is positive semidefinite.

Clearly $\hmx$ is a nonnegative matrix. In fact:

\begin{lemma}
    If $X$ is connected, all entries of $\hmx$ are positive.
\end{lemma}

\proof
Each Schur square $E_r^{\circ2}$ is nonnegative and if $(\hmx)_{u,v}=0$
then $(E_r^{\circ2})_{u,v}=0$. However this implies that $(E_r)_{u,v}=0$
for all $r$ and hence any linear combination of the idempotents $E_r$ has $uv$-entry
zero. Since this implies that $(A^k)_{u,v}=0$ for all $k$ we conclude
that $X$ is not connected.\qed

In \cite[Lemma 16.2]{Godsil2011a} we show that average mixing is never uniform,
that is, $\hmx$ cannot be a scalar multiple of the all-ones matrix $J$.

We say that a matrix $M$ is a \textsl{contraction} if $x^*Mx \le x^*x$ for all
complex vectors $x$.

\begin{lemma}
    If $\hmx$ is the average mixing matrix of the graph $X$ then all eigenvalues
    of $\hmx$ lie in the interval $[0,1]$.
\end{lemma}

\proof
Since $H(t)\otimes H(-t)$ is unitary and since $H(t)\circ H(-t)$ is a principal submatrix
of $H(t)\otimes H(-t)$, it follows that $H(t)\circ H(-t)$ is a contraction. Since it is
symmetric and real, its eigenvalues must lie in the interval $[0,1]$, and therefore
this holds true for $\hmx$ too.\qed

Since $\hmx$ is the average of doubly stochastic matrices, it is doubly stochastic
and its largest eigenvalue is 1. We can also see this without appealing to the 
averaging. Notice that
\[
    ((E_r\circ E_r)\one)_u = ((e_u^TE_r)\circ(e_u^TE_r))\one = \ip{e_u^TE_r}{e_u^TE_r}
        =e_u^TE_r^2e_u =e_u^TE_re_u =(E_r)_{u,u}
\]
and, since $\sum_r E_r=I$, it follows that $\hmx\one=\one$.

\begin{lemma}
    The average mixing matrix of a graph is rational.
\end{lemma}

\proof
Let $\phi(X,x)$ be the characteristic polynomial of $X$, and let $\fld$ be a splitting 
field for $\phi(X,x)$. We use the fact that an element of $\fld$ which is fixed
by all field automorphisms of $\fld$ must be rational. If $\sg$ is an automorphism
of $\fld$, then
\[
    A = A^\sg = \sum_r \th_r^\sg E_r^\sg.
\]
Since $\th_r^\sg$ must be an eigenvalue of $A$ and since the spectral decomposition
of $A$ is unique, it follows that $E_r^\sg$ is one of the idempotents in the spectral
decomposition of $A$. Therefore the set of idempotents is closed under field
automorphisms, and so must the $\{E_r^{\circ2}\}_r$. Consequently
\[
    \hmx^\sg = \hmx
\]
for all $\sg$ and therefore $\hmx$ is rational.\qed

Note that this lemma holds whether we use the Laplacian or the adjacency 
matrix---all we need is that $A$ be symmetric with integer entries.

We use $L(X)$ to denote the Laplacian matrix $X$. If $\De$ is the diagonal matrix 
whose $i$-th diagonal entry is the valency of the $i$-vertex of $X$, then
\[
    L(X) =\De - A.
\]
The transition matrix of $X$ relative to $L(X)$ is
\[
    H_L(t) := \exp(it(\De-A)).
\]
When $X$ is regular, questions about $H_L$ reduce immediately to questions about
$H_X$, but in general there is no simple relation between the two cases.

\begin{lemma}
    If $X$ is regular then $X$ and its complement $\comp{X}$ have the same 
    average mixing matrix. For any graph $X$, the average mixing matrix relative 
    to the Laplacian of $X$ is equal to the average mixing matrix relative to 
    the Laplacian of $\comp{X}$.
\end{lemma}

\proof
If $X$ is regular then the idempotents in the spectral decomposition of its adjacency
matrix are the idempotents in the spectral decomposition of the adjacency matrix of 
$\comp{X}$. For any graph $X$ on $n$ vertices
\[
    L(\comp{X}) = L(K_n) -L(X);
\]
since $L(K_n)=nI-J$ and since $L(X)$ commutes with $J$, 
the idempotents in the spectral decomposition of its 
Laplacian are the idempotents in the spectral decomposition of the Laplacian of 
$\comp{X}$.\qed

\section{Integrality}

In investigating the relation between the structure of $\hmx$ and the graph $X$,
it can be convenient to scale $\hmx$ so that it entries are integers. For this we need
to know a common multiple of the denominators of its entries.

\begin{lemma}
    If $D$ is the discriminant of the minimal polynomial of $A$, then $D^2\hmx$
    is an integer matrix.
\end{lemma}

\proof
Let $\seq\th1m$ be the distinct eigenvalues of $A$. Define polynomials
$\ell_r(t)$ by
\[
    \ell_r(t) := \prod_{s\ne r}(t-\th_s).
\]
We note that $\ell_r(\th_r)=\psi'(\th_r)$ and
\[
    \frac{\ell_r(\th_s)}{\psi'(\th_r)} = \de_{r,s}.
\]
Now
\[
    E_r = \frac1{\psi'(\th_r)} \ell_r(A).
\]
The discriminant $D$ of $\psi$ is equal (up to sign) to
\[
    \prod_{r=1}^m \psi'(\th_r);
\]
since the entries of $\ell_r(A)$ are algebraic integers we conclude that the entries
of $D^2E_r^{\circ2}$ are algebraic integers and therefore the entries of $D^2\hmx$
are algebraic integers. Since $\hmx$ is rational, the lemma follows.\qed

We have no reason to believe this lemma is optimal. If the eigenvalues of $X$ 
are all simple, we can do better.

\begin{theorem}
    Let $X$ be a graph with all eigenvalues simple and let $D$ be the discriminant
    of its characteristic polynomial. Then $D\hmx$ is an integer matrix. 
\end{theorem}

\proof
We have
\[
    (xI-A)^{-1} = \sum_r \frac1{x-\th_r}E_r
\]
and since $(I-tA)^{-1}$ is the walk generating function of $X$, it follows
from \cite[Corollary~4.1.3]{cgblue} that
\[
    (E_r)_{u,v} =\lim_{x\to\th_r} 
        \frac{
        (x-\th_r)
        \bigl(\phi(X\diff u,x)\phi(X\diff v,x) -\phi(X\diff uv,x)\phi(X,x)\bigr)^{1/2}}
        {\phi(X,x)}
\]
and since
\[
    \lim_{x\to\th_r}\frac{\phi(X,x)}{x-\th_r} = \phi'(X,\th_r)
\]
we conclude that if $\th_r$ is simple
\[
    \bigl((E_r)_{u,v}\bigr)^2 
        = \frac{\phi(X\diff u,\th_r)\,\phi(X\diff v,\th_r)}{\phi'(X,\th_r)^2}.
\]
If $B$ is the $n\times n$ matrix with $ur$-entry $\phi(X\diff u,\th_r)$ and $\De$
is the $n\times n$ diagonal matrix with $r$-th diagonal entry $\phi'(X,\th_r)$, 
it follows that
\[
    \hmx = B\De^{-2}B^T.
\]

Assume $n=|V(X)|$ and let $\seq\th1n$ be the eigenvalues of $X$. Let $V$ be the
$n\times n$ Vandermonde matrix with $ij$-entry $\th_j^{i-1}$. Let $\phi$
be the characteristic polynomial of $X$. The discriminant
of $\phi$ is equal to the product of the entries of $\De$, we denote it by $D$.

Let $C$ be the $n\times n$ matrix whose $ur$-entry is the coefficient of 
$x^{r-1}$ in $\phi(X\diff u,x)$. Then $CV=B$ and
\[
    \hmx = CV\De^{-2}V^TC^T.
\]
Define polynomials $\ell_s(x)$ by
\[
    \ell_s(x) = \prod_{r\ne s} (x-\th_r)
\]
and let $L$ be the $n\times n$ matrix with $sj$-entry equal to the coefficient of
$x^{j-1}$ in $\ell_s(x)$. Note that
\[
    \ell_s(\th_r) = \de_{r,s}\phi'(\th_s)
\]
and therefore
\[
    LV = \De.
\]
Since
\[
    \De = \De^T = V^TL^T
\]
it follows that
\[
    \De^{-2} = V^{-1}L^{-1}L^{-T}V^{-T}
\]
and so
\[
    V\De^{-2}V^T = VV^{-1}L^{-1}L^{-T}V^{-T}V^T = (L^TL)^{-1}.
\]

(We're almost done.) The entries of $L$ are algebraic integers. As
\[
    \det(VV^T) = \det(V)^2 = D
\]
and as $LD=\De$, we see that $\det(L^TL)=D$. Therefore the entries of
\[
    D(L^TL)^{-1}
\]
are algebraic integers. So the entries of $DV\De^{-2}V^T$ are algebraic integers
and, since the entries of $C$ are integers, the entries of $D\hmx$ are algebraic
integers. But the entries of $D\hmx$ are rational and therefore they are all 
integers.\qed

It is at least plausible that if $D$ is the discriminant of the minimal polynomial
of $X$, then $D\hmx$ is an integer matrix.

\section{Average Mixing on Paths}

We need the following trigonometric identity.

\begin{lemma}
    \label{lem:cossum}
    \[
        2\sum_{r=0}^n \cos(r\th) = \frac{\sin((n+\frac12)\th)}{\sin(\frac12\th)}+1.
    \]
\end{lemma}

\proof
If $q:=e^{i\th}$ then
\begin{align*}
    2\sum_{r=0}^n \cos(r\th) &= \sm r0n(q^r+q^{-r})\\
        &= \frac{q^{n+1}-1}{q-1} + \frac{q^{-n-1}-1}{q^{-1}-1}\\
        &= \frac{q^{n+1}-1}{q-1} +\frac{q^{-n}-q}{1-q}\\
        &= \frac{q^{n+1}-q^{-n}}{q-1} +1\\
        &= \frac{q^{n+1/2}-q^{-1/2-n}}{q^{1/2}-q^{-1/2}}+1\\
        &= \frac{\sin((n+\frac12)\th)}{\sin(\frac12\th)}+1.\tag*{\sqr53}
\end{align*}

We will also use the following explicit expression for the idempotents in the spectral
decomposition of $A(P_n)$. This result is standard, but we offer a sketch of the proof.

\begin{lemma}
    \label{lem:pnids}
    The idempotents $\seq E1n$ in the spectral decomposition of $P_n$ are given by
    \[
        (E_r)_{j,k} = \frac2{n+1} \sin\left(\frac{jr\pi}{n+1}\right)
            \sin\left(\frac{kr\pi}{n+1}\right).
    \]
\end{lemma}

\proof
If $A=A(P_n)$ and $e_n$ is the $n$ vector in the standard basis of $\re^n$, then
\[
    A\pmat{\sin(\be)\\ \sin(2\be)\\ \vdots\\ \sin(n\be)} 
        = \pmat{\sin(2\be)\\ \sin(\be)+\sin(3\be)\\ \vdots\\ \sin((n-1)\be)}
        = 2\cos(\be)\pmat{\sin(\be)\\ \sin(2\be)\\ \vdots\\ \sin(n\be)}
            -\sin((n+1)\be)e_n
\]
So if $\sin((n1\be))=0$ then
\[
    z(\be) := \pmat{\sin(\be)\\ \sin(2\be)\\ \vdots\\ \sin(n\be)}
\]
is an eigenvector for $A$ with eigenvalue $2\cos(\be)$. Letting $\be$ vary over
the values
\[
    \frac{2\pi r}{n+1},\quad r=1,\ldots,n
\]
we obtain $n$ distinct eigenvalues. Therefore each eigenvalue of $P_n$ is simple
and the projection onto the eigenspace spanned by $z(\be)$ is
\[
    \frac1{z(\be)^Tz(\be)} z(\be)z(\be)^T).
\]
We can compute the value of the inner product $z(\be)^Tz(\be)$ using \lref{cossum},
and this yields the stated expression for $E_r$.\qed

\begin{lemma}
    If $\seq E1n$ are the idempotents for $P_n$, then the average mixing 
    matrix of $P_n$ is
    \[
        \sum_r E_r\circ E_r = \frac1{2n+2}(2J+I+T).
    \]
\end{lemma}

\proof
We use \lref{pnids}:
\[
    (E_r\circ E_r)_{j,k} = \frac4{(n+1)^2} \sin^2\Bigl(\frac{jr\pi}{n+1}\Bigr) 
        \sin^2\Bigl(\frac{kr\pi}{n+1}\Bigr)
\]
which implies that
\[
    \frac{(n+1)^2}4 (E_r\circ E_r)_{j,k} 
        = \frac14\left(1-\cos\Bigl(\frac{2jr\pi}{n+1}\Bigr)\right)
        \left(1-\cos\Bigl(\frac{2kr\pi}{n+1}\Bigr)\right).
\]
Now
\begin{multline*}
    \left(1-\cos\Bigl(\frac{2jr\pi}{n+1}\Bigr)\right)
        \left(1-\cos\Bigl(\frac{2kr\pi}{n+1}\Bigr)\right) =\\
    1-\cos\Bigl(\frac{2jr\pi}{n+1}\Bigr) -\cos\Bigl(\frac{2kr\pi}{n+1}\Bigr)
        +\frac12\cos\Bigl(\frac{2(j+k)r\pi}{n+1}\Bigr)
        +\frac12\cos\Bigl(\frac{2(j-k)r\pi}{n+1}\Bigr)
\end{multline*}
We need to sum each of the five terms on the right from 1 to $n$. From
\lref{cossum} it follows that
\[
    \sm r1n \cos\Bigl(\frac{2\ell r\pi}{n+1}\Bigr) 
        = \frac12\left(-1+\frac{\sin\Bigl(\frac{(2n+1)\ell\pi}{n+1}\Bigr)}
            {\sin\Bigl(\frac{\ell\pi}{n+1}\Bigr)}\right)
        = \frac12\left(-1+\frac{\sin\Bigl(2\ell\pi-\frac{\ell\pi}{n+1}\Bigr)}
            {\sin\Bigl(\frac{\ell\pi}{n+1}\Bigr)}\right)
        = -1.
\]
Consequently
\[
    \sm r1n (n+1)^2 (E_r\circ E_r)_{j,k} = \begin{cases}
        3(n+1)/2,& j=k;\\
        3(n+1)/2,& j+k=n+1\\
        n+1,& \mathrm{otherwise}
    \end{cases}
\]
and this completes the proof.\qed

\section{Path Laplacians}

Let $\De$ denote the diagonal matrix with $i$-th diagonal entry equal to the
valency of the $i$-th vertex of $X$. Then the Laplacian of $X$ is $\De-A$.
Let $D$ be the $n\times(n-1)$ matrix with
\[
    D_{i,i}=1,\quad D_{i,i-1}=-1
\]
and all other entries zero. Then $D$ is the incidence matrix of an orientation
of $P_n$ and
\[
    DD^T =\De- A(P_n),\qquad D^TD = 2I-A(P_{n-1}).
\]
If $\seq\th1{n-1}$ are the eigenvalues of $P_{n-1}$, this shows that the non-zero
eigenvalues of $\De-A(P_n)$ are the numbers $2-\th_r$ for $r=1,\ldots,n-1$.
We can also use this to determine the idempotents in the spectral decomposition
of $\De-A(P_n)$. The lemma follows.\qed

\begin{lemma}
   If $\seq E1{n-1}$ are the idempotents in the spectral decomposition of $P_{n-1}$,
   then the idempotents of $\De-A(P_n)$ are $n^{-1}J$ and
   \[
        \frac1{2-\th_r}DE_rD^T,\ldots, r=1,\ldots,n-1.
   \]
\end{lemma}

\proof
Since
\[
    DE_r D^TDE_sD^T = DE_r(2I-A(P_{n-1})E_s D^T)
\]
and since $E_rE_s=0$ if $r\ne s$ and $(2I-A(P_{n-1})E_s =(2-\th_s)E_s$, it follows
that
\[
    DE_r D^TDE_sD^T =\de_{r,s}(2-\th_r)DE_rD^T.
\]
Therefore $(2-\th_r)^{-1}DE_rD^T$ is an idempotent. Next
\[
    (\De-A(P_n))DE_rD^T = DD^T DE_rD^T = D(2I-A(P_{n-1}))E_rD^T
        = (2-\th_r) DE_rD^T
\]
and therefore $(2-\th_r)^{-1}DE_rD^T$ represents orthogonal projection onto
an eigenspace of $A(P_n)$. The lemma follows.\qed

\begin{lemma}
    If $\seq E1{n-1}$ are the idempotents in the spectral decomposition of $P_{n-1}$,
    then
    \[
        (2-\th_r)^{-1}(DE_rD^T)_{j,k} 
            = \frac{2}{n} 
            \cos\left(\frac{(2j-1)r\pi}{2n}\right) 
            \cos\left(\frac{(2k-1)r\pi}{2n}\right), \qquad 1\le j,k\le n.
    \]
\end{lemma}

\proof
From \lref{pnids}, we have
\[
    (E_r)_{j,k} = \frac2{n} \sin\left(\frac{jr\pi}{n}\right)
        \sin\left(\frac{kr\pi}{n}\right), \qquad 1\le j,k\le n-1.
\]
Let $\al=r\pi/n$ and let $\sg$ denote the column vector of length $n-1$ where 
$\sg_j=\sin(j\al)$. Then
\[
    DE_rD^T = \frac2{n} D\sg (D\sg)^T
\]
and 
\[
    D\sg =\pmat{
        \sin(\al)\\ 
        \sin(2\al)-\sin(\al)\\
        \vdots\\
        \sin((n-1)\al)-\sin((n-2)\al)\\
        -\sin((n-1)\al)}
        = 2\sin(\al/2)\pmat{
            \cos(\al/2)\\
            \cos(3\al/2)\\
            \vdots\\
            \cos((2n-3)\al/2)\\
            \cos((2n-1)\al/2)},
\]
where in computing the last entry we have used the fact that $n\al=r\pi$,
whence $\sin(n\al)=0$ and
\[
    -\sin(n-1)\al = \sin(n\al) - \sin(n-1)\al.
\]

Finally for $P_{n-1}$ we have
\[
    2-\th_r = 2-2\cos\Bigl(\frac{r\pi}{n}\Bigr) = 4\sin^2\Bigl(\frac{r\pi}{2n}\Bigr)
\]

\begin{theorem}
    The average mixing matrix for the continuous quantum walk using the
    Laplacian matrix of the path $P_n$ is
    \[
        \frac1{n^2}\left((n-1)J+\frac{n}{2}(I+T)\right).
    \]
\end{theorem}

\proof
Set 
\[
    F_r = (2-\th_r)^{-1}DE_rD^T.
\]
Then
\begin{align*}
    (F_r^{\circ2})_{j,k} &= \frac{4}{n^2}\cos^2\left(\frac{(2j-1)r\pi}{2n}\right) 
        \cos^2\left(\frac{(2k-1)r\pi}{2n}\right)\\
    &= \frac{1}{n^2} \left(1+\cos\left(\frac{(2j-1)r\pi}{n}\right)\right)
            \left(1+\cos\left(\frac{(2k-1)r\pi}{n}\right)\right)
\end{align*}
and
\begin{align*}
    \left(1+\cos\Bigl(\frac{(2j-1)r\pi}{n}\Bigr)\right)
          & \left(1 +\cos\Bigl(\frac{(2k-1)r\pi}{n}\Bigr)\right)\\[2pt]
      = 1 &+ \cos\Bigl(\frac{(2j-1)r\pi}{n}\Bigr) 
            + \cos\Bigl(\frac{(2k-1)r\pi}{n}\Bigr)\\[2pt]
            &+ \frac12\cos\Bigl(\frac{(2j+2k-2)r\pi}{n}\Bigr)
            + \frac12\cos\Bigl(\frac{(2j-2k)r\pi}{n}\Bigr).
\end{align*}

From \lref{cossum} we have
\begin{align*}
    2\sm{r}1{n-1}\cos\Bigl(\frac{r\ell\pi}{n}\Bigr) 
    &= -1 + \frac{\sin\bigl((n-\frac12)\frac{\ell\pi}{n}\bigr)}
            {\sin\bigl(\frac{\ell\pi}{2n}\bigr)}\\[2pt]
    &= -1+\frac{\sin\bigl(\ell\pi-\frac{\ell\pi}{2n}\bigr)}
        {\sin\bigl(\frac{\ell\pi}{2n}\bigr)}\\[2pt]
    &= -1 + \frac{-\cos(\ell\pi)\sin\Bigl(\frac{\ell\pi}{2n}\Bigr)}
        {\sin\Bigl(\frac{\ell\pi}{2n}\Bigr)}\\[2pt]
    &= -((-1)^\ell+1).
\end{align*}
It is now easy to derive the stated formula for the average mixing matrix.\qed

We note that $2I-L(P_n)$ can be viewed as the adjacency matrix of a path on $n$ 
vertices with a loop of weight one on each end-vertex. Examples show that if we add 
loops with weight other than 0 or 1, the average mixing matrix is not a linear 
combination of $I$, $J$ and $T$. Thus if we add loops of weight 2 to the end-vertices 
of $P_6$, the average mixing matrix is:
\[
   \frac1{2*9*107} 
   \left(\begin{array}{rrrrrr}
   599 & 218 & 146 & 146 & 218 & 599 \\
   218 & 455 & 290 & 290 & 455 & 218 \\
   146 & 290 & 527 & 527 & 290 & 146 \\
   146 & 290 & 527 & 527 & 290 & 146 \\
   218 & 455 & 290 & 290 & 455 & 218 \\
   599 & 218 & 146 & 146 & 218 & 599
   \end{array}\right). 
\]
(Here the discriminant of the characteristic polynomial is $2^6  3^5  107$.)

\section{Cycles}

We determine the average mixing matrices for cycles.

Let $P$ be the permutation
matrix corresponding to a cycle of length $n$ and let $\zeta$ be a primitive
$n$-th root of unity. Define matrices $\seq F0{n-1}$ by
\[
    (F_r)_{i,j} = \frac1n\zeta^{r(i-j)}.
\]
(Thus the rows and columns of these matrices are indexed by $\{0,\ldots,n-1$.)
Then $P$ is a normal matrix and has the spectral decomposition
\[
    P = \sum_{r=0}^{n-1} \zeta^r F_r.
\]
We also note that
\[
    F_r\circ F_s = \frac1{n} F_{r+s}
\]
where the subscripts are viewed as elements of $\ints_n$. The adjacency matrix
of the cycle on $n$ vertices is $P+P^T$. Define $E_0$ to be $F_0$ and, 
if $0<r<n/2$, we set
\[
    E_r = F_r + F_{n-r}.
\]
Further $E_0:=F_0$ and, if $n$ is even then $E_{n/2}:=F_{n/2}$.
Then if $n$ is odd, 
\[
    \seq E0{(n-1)/2}
\]
are the idempotents in the spectral decomposition 
of $A(C_n)$ and the corresponding eigenvalues are
\[
    \th_r  = \zeta^r+\zeta^{-r}, \qquad r=0,\ldots,\frac{n-1}2;
\]
if $n$ is even we have the additional idempotent $E_{n/2}$ with eigenvalue $\zeta^{n/2}=-1$.

\begin{lemma}
    If $n$ is odd then the average mixing matrix of the cycle $C_n$ is 
    \[
        \frac{n-1}{n^2}J + \frac1{n}I,
    \]
    if $n$ is even it is 
    \[
        \frac{n-2}{n^2}J +\frac{1}{n}(I+P^{n/2}).
    \]
\end{lemma}

\proof
Assume first that $n=2m+1$. Then the average mixing matrix is
\begin{align*}
    \sum_{r=0}^m E_r^{\circ2} &= F_0^{\circ2} +\sum_{r=1}^m (F_r+F_{-r})^{\circ2}\\[2pt]
        &= \frac{1}{n}F_0 + \frac1{n}\sum_{r=1}^m (F_{2r}+F_{-2r}+2F_0)\\[2pt]
        &= F_0 + \frac1{n}\sum_{r=1}^{n-1} F_r\\[2pt]
        &= \frac1n I + \frac{n-1}{n}F_0.
\end{align*}

Now suppose $n=2m$. Then the average mixing matrix is
\begin{align*}
    E_m^{\circ2}+\sum_{r=0}^{m-1}E_r^{\circ2} &= F_0^{\circ2}+F_m^{\circ2}
            +\sum_{r=1}^{m-1}(F_r+F_{-r})^{\circ2}\\[2pt]
        &= \frac{1}{n}F_0 +\frac{1}{n}F_0 +\frac{1}{n}  
            \sum_{r=1}^{m-1}(F_{2r}+F_{-2r}+2F_0)\\[2pt]
        &= F_0 +\frac1{n}\sum_{r=1}^{m-1}(F_{2r}+F_{-2r})\\[2pt]
        &= \frac{n-2}{n}F_0 +\frac2{n}\sum_{r=0}^{m-1}F_{2r}.
\end{align*}
Since
\[
    P^m F_s = (\zeta^m)^s F_s= (-1)^s F_s
\]
we see that $P^m$ has spectral decomposition
\[
    P^m =\sum_{s=0}^{n-1} (-1)^s F_s
\]
and consequently
\[
    \frac12(I+P^m) = \sum_{r=0}^{m-1}F_{2r}.
\]
Our stated formula for the average mixing matrix follows.\qed

In general the average mixing matrix for a circulant can be more complex in 
structure than the average mixing matrix of a cycle. We will see in the next
section that graphs in pseudocyclic schemes provide the right generalization
of cycles or, at least, of odd cycles.

It is easy for different circulants of order $n$ to have the same spectral 
idempotents, any two such circulants necessarily have the same average mixing matrix.

\section{Pseudocyclic Schemes}

An association scheme $\cA$ with $d$ classes on $v$ vertices is a set of
$01$-matrices $\seq A0d$ such that
\begin{enumerate}[(a)]
    \item 
    $A_0=I$ and $\sum_i A_i=J$.
    \item
    For each $i$ we have $A_i^T\in\cA$.
    \item
    $A_iA_j=A_jA_i$.
    \item
    For all $i$ and $j$ the product $A_iA_j$ lies in the real span of $\cA$.
\end{enumerate}
Convenient references are \cite{bcn, cgblue}. 
We often view the matrices $A_i$ as the adjacency matrices of directed graphs on $v$ 
vertices, we call these the \textsl{graphs} of the scheme and note that the
axioms are often stated in terms of the graphs rather than their adjacency matrices.
An association scheme is \textsl{symmetric} if each matrix in it is symmetric. These
axioms imply that the span of $\cA$ is a commutative matrix algebra that is
closed under the Schur product. It is called the \textsl{Bose-Mesner} algebra of the
scheme. The matrices in $\cA$ form a basis consisting of Schur idempotents. There is
also a basis of matrix idempotents $\seq E0d$, which satisfy
\begin{enumerate}[(a)]
    \item
    $E_0 = v^{-1}J$ and $\sum E_j=I$.
    \item
    For each $j$ there is an index $j'$ such that $\comp{E}_j= E_{j'}$.
    \item
    For all $i$ and $j$ the Schur product $E_i\circ E_j$ lies in the Bose-Mesner algebra.
\end{enumerate}
The matrix idempotents can be viewed as providing (an analog of) the spectral 
decomposition.

Since the matrices $A_i$ commute it follows that $A_i$ and $J$ commute and hence
there are constants $v_i$ such that $AJ=v_i J$, these are the \textsl{valencies}
of the scheme. The column space of $E_j$ is an eigenspace for all matrices in the
Bose-Mesner algebra, its dimension is denoted by $m_j$ and it is referred to
as a \textsl{multiplicity} of the scheme. Note that $v_0=m_0=1$.
A scheme is \textsl{pseudocyclic} if $\seq m1d$ are equal (in which case their common
value is $(v-1)/d$). If a scheme is pseudocyclic then $\seq v1d$ are necessarily
all equal to $(v-1)/d$. For details see Brouwer, Cohen and Neumaier\cite[\S2.2B]{bcn}. 

We note one class of examples, the \textsl{cyclotomic schemes}.
Assume $q$ is a prime power and $d$ divides $q-1$. Let $\fld$ be the finite 
field of order 
$q$ and let $S$ be the subgroup of the multiplicative group of $\fld$ generated
by the non-zero $d$-th powers. Thus $|S|=(q-1)/d$. Let $\seq S1d$ denote the cosets of
$S$ in $\fld^*$, with $S_1=S$. Now we define the matrices of an association scheme
with $d$ classes and with vertex set $\fld$ by setting $A_0=I$ and
\[
    (A_i)_{x,y} = 1,\qquad (i=1,\ldots,d)
\]
if and only if $y-x$ is in $S_i$. These matrices form the cyclotomic scheme with $d$
classes on $\fld$. It is symmetric if and only if $-1\in S$. The directed graphs 
$\seq X1d$ such that $A_i=A(X_i)$ are all isomorphic. The most well known
case is when $d=2$ and $q\cong1$ modulo four, in which case the two graphs we 
have constructed are the \textsl{Paley graphs}.

We note that there are pseudocyclic schemes which are not cyclotomic, and that there
are pseudocyclic schemes with two classes where the two graphs are asymmetric, that is,
their only automorphism is the identity.

\section{Average Mixing on Pseudocyclic Graphs}

\begin{theorem}
    Suppose $X$ is a graph in a $d$-class pseudocyclic scheme on $n$ vertices consisting
	of graphs of valency $m=(n-1)/d$. Then the average mixing matrix of $X$ is
    \[
        \frac{n-m+1}{n^2}J + \frac{m-1}{n}I.
    \]
\end{theorem}

\proof
Koppinen \cite{koppinen} proved that for an association scheme with $d$ classes we have
\[
    \sm{i}0d \frac1{nv_i}A_i\otimes A_i^T = \sm{j}0d \frac1{m_j}E_j\otimes E_j.
\]
If the scheme is symmetric then $A_i\circ A_i=A_i$ and, since $M\circ N$ is a submatrix 
of $M\otimes N$, Koppinen's identity yields that
\[
    \sm{i}0d \frac1{nv_i}A_i = \sm{j}0d \frac1{m_j}E_j\circ E_j.
\]
For any scheme, we have $m_0=v_0=1$ and for a pseudocyclic scheme 
\[
    m_i=v_i=\frac{n-1}{d},\quad i=1,\ldots,d.
\] 
As $\sum_i A_i=J$ and $\sum_j E_j=I$, if $m=(n-1)/d$ we find that
\[
    \frac1n\Bigl(I +\frac1m(J-I)\Bigr) 
        = \frac1{n^2}J+\frac1m\sm{j}1d E_j^{\circ2}\Bigr)
\]
and hence
\[
    \frac{1}{m}\sm{j}1d E_j^{\circ2} 
        = \left(\frac1{nm}-\frac1{n^2}\right)J + \frac{m-1}{nm}I
        =\frac{n^2-nm}{n^3m}J + \frac{m-1}{nm}I.\qed
\]

A graph in a pseudocyclic scheme with two classes is known as a 
\textsl{conference graph}, and there are examples whose automorphism group is trivial.

\section{Cospectral Vertices}

If $u\in V(X)$ then $(A^k)_{u,u}$ is equal to the number of closed walks on $X$
of length $k$ that start (and finish) at $u$. The following result summarizes some 
standard facts. (The equivalence of (a) and (c) follows from Equation~(1) on 
page~52 of \cite{cgblue} and (a) and (b) from Equation~(2) on page~9 of the same source.)

\begin{theorem}\label{thm:covx}
    Let $X$ be a graph with vertices $u$ and $v$, and let $\seq E1m$ be the 
    idempotents in the spectral decomposition of $A(X)$.
    Then the following are equivalent:
    \begin{enumerate}[(a)]
        \item 
        $(A^k)_{u,u}=(A^k)_{v,v}$ for all $k\ge0$.
        \item
        $\|E_re_u\|^2=\|E_re_v\|^2$.
        \item
        The graphs $X\diff u$ and $X\diff v$ are cospectral.\qed
    \end{enumerate}
\end{theorem}

Note that
\[
    \|E_re_u\|^2 = e_u^TE_r^TE_re_u = e_u^T E_r^2e_u =e_u^TE_re_u =(E_r)_{u,u}.
\]
We will say that vertices $u$ and $v$ of $X$ are \textsl{cospectral} if any of
the conditions in this theorem hold.

We recall that the graph $X$ admits \textsl{perfect state transfer} from $u$
to $v$ at time $\tau$ if there is a scalar $\ga$ such that
\[
    H_X(\tau)e_u  = \ga e_v
\]
where $|\ga|=1$ since $\|e_u\| =\|H_X(\tau)e_u\|$. As $H_X(t)$ is a polynomial in $A$,
we have
\[
    E_rH_X(t) = H_X(t)E_r =\exp(i\th_r t)E_r
\]
and therefore
\[
    \ga E_re_v = H_X(\tau) E_re_u = \exp(i\th_r\tau) E_re_v.
\]
Since $E_re_v$ and $E_re_u$ are real vectors and since $|\ga|=|\exp(i\th_r\tau)|=1$,
we conclude that
\[
    E_r e_u = \pm E_r e_v.
\]
(This fact is well known and we include the proof for convenience only.)

Clearly if $E_r e_u = \pm E_r e_v$ then
\[
    \|E_re_u\|^2=\|E_re_v\|^2
\]
and so $u$ and $v$ are cospectral. We will say that $u$ and $v$ are
\textsl{strongly cospectral} if $E_r e_u = \pm E_r e_v$.

\begin{lemma}
    If $u$ and $v$ are vertices in $X$ and the eigenvalues of $X$ are simple, then
    $u$ and $v$ are strongly cospectral if and only if they are cospectral.
\end{lemma}

\proof
We only need to consider the case where $u$ and $v$ are cospectral. The vectors
$E_re_u$ and $E_re_v$ are both eigenvectors of $X$ with eigenvalue $\th_r$.
Since $\th_r$ is simple this means that that $E_re_v$ is a scalar multiple of
$E_re_u$. Since these two vectors have the same length, the lemma follows.\qed

\begin{theorem}
    Let $\hmx$ be the average mixing matrix of the graph $X$. Then vertices
    $u$ and $v$ are strongly cospectral if and only if $\hmx(e_u-e_v)=0$.
\end{theorem}

\proof
Suppose $N\psd0$ and $N(e_1-e_2)=0$. We may assume that the leading $2\times2$
submatrix of $N$ is
\[
    \pmat{a&b\\ b&d}
\]
and therefore
\[
    0 = (e_1-e_2)^TN(e_1-e_2) = a+d-2b.
\]
Hence $b=(a+d)/2$. Since $N\psd0$ we have $ad-b^2\ge0$ and thus
\[
    0 \le 4ad -4b^2 = 4ad - (a+d)^2 = -(a-d)^2,
\]
whence $a=d=b$.

If $\hmx(e_u-e_v)=0$ then
\[
    0 = (e_u-e_v)^T\hmx(e_u-e_v) = \sum_r (e_u-e_v)^T E_r^{\circ 2}(e_u-e_v)
\]
and as each summand $E_r^{\circ2}$ in $\hmx$ is positive semidefinite, we have
\[
    E_r^{\circ 2}(e_u-e_v) = 0
\]
for all $r$. Therefore, for all $r$,
\[
    ((E_r)_{u,u})^2 = ((E_r)_{u,v})^2 =((E_r)_{v,v})^2
\]
Since
\[
    (E_r)_{u,v} = \ip{E_re_u}{E_re_v}
\]
it follows by Cauchy-Schwarz that $E_re_u=\pm E_re_v$.\qed

A graph is \textsl{walk regular} if all its vertices are cospectral. There
are many classes of examples, and it is not necessary that such graphs
be vertex transitive. We note the following though:

\begin{lemma}
    If all vertices in $X$ are strongly cospectral, then $X$ is $K_1$ or $K_2$.
\end{lemma}

\proof
If all vertices are strongly cospectral, then all rows of $\hmx$ are equal,
and therefore it is a scalar multiple of $J$. By \cite[Lemma 16.2]{Godsil2011a} this
implies that $|V(X)|\le2$.\qed

One consequence of the results here is that, if the rows of the average mixing matrix
of $X$ are distinct, there is no perfect state transfer on $X$. 
Kay \cite[Section~D]{Kay2011} has shown that if we have perfect state transfer in $X$ 
from $u$ to $v$ and from $u$ to $w$, then $v=w$.

\section{Discrete Walks}

If $U$ is a unitary matrix then its powers describe a \textsl{discrete
quantum walk}. The corresponding average mixing matrix is
\[
	\lim_{N\to\infty} \sum_{n=0}^{N-1} U^n\circ U^{-n}.
\]
Since unitary matrices are normal, they have a spectral decomposition which
we can write as
\[
	U = \sum_r e^{i\th_r}E_r
\]
where $\th_r$ is real and $E_r$ Hermitian for all $r$. The unitary matrices
considered in quantum computing are constructed from an underlying regular
graph.

\begin{theorem}
	If $\seq E1m$ are the idempotents in the spectral decomposition of the unitary 
	matrix $U$, the average mixing matrix $\widehat U$ of $U$ is $\sum_r E_r^{\circ2}$.
	If $U$ is rational then $\hat U$ is rational.
\end{theorem}

\proof
We have
\[
    U^n\circ U^{-n} = \sum_r E_r^{\circ2} +\sum_{r\ne s}e^{ni(\th_r-\th_s)} E_r\circ E_s.
\]
Therefore
\[
    \frac1{N}\sum_{n=0}^{N-1} U^n\circ U^{-n} = \sum_r E_r^{\circ2}
        + \frac1{N}\left(\sum_{n=0}^{N-1} e^{ni(\th_r-\th_s)}\right)E_r\circ E_s
\]
and since
\[
    \sum_{n=0}^{N-1} e^{ni(\th_r-\th_s)} = \frac{e^{ni(\th_r-\th_s)}-1}{e^{i(\th_r-\th_s)}-1}
\]
The right side is bounded in absolute value by
\[
    \frac2{\left|e^{i(\th_r-\th_s)}-1\right|}
\]
which is independent of \(N\). Hence
\[
    \frac1{N}\sum_{n=0}^{N-1} U^n\circ U^{-n} \to \sum_r E_r^{\circ2}
\]
as $N$ tends to infinity.\qed

Since
\[
	I - U^n\circ U^{-n} = \sum_{r\ne s}(1-e^{ni(\th_r-\th_s)}) E_r\circ E_s
		=2\sum_{r\ne s}(1-\cos(n(\th_r-\th_s))) E_r\circ E_s,
\]
we see that $I-U^n\circ U^{-t}$ is positive semidefinite and hence $I-\widehat{U}$
is positive semidefinite.

If $U$ is based on a regular graph then $\widehat{U}$ is rational.

\section{Questions}

We list three questions raised by our work. 

\begin{enumerate}
    \item
    Is it true that if $D$ is the discriminant of the minimal polynomial of $X$,
    then $D\hmx$ is an integer matrix? 
    \item
    Is there a useful algorithm for computing $\hmx$ that works over the rationals? 
    \item
    Which graph have the property that their average mixing matrix is a linear combination 
    of $I$ and $J$? (Jamie Smith, in unpublished work, has found a strongly regular graph 
    with this property that is not pseudocyclic.)
\end{enumerate}

\section*{Acknowledgements}

I have had some very useful discussions about average mixing with Christino Tamon.
All computations were carried out using the package sage \cite{sage}. 


\end{document}